\documentclass[12pt]{article}
\usepackage{amssymb}
\font\smallit=cmti10
\font\smalltt=cmtt10
\font\smallrm=cmr9

\usepackage{epsfig}

\newenvironment{packed_enumerate}{
\setlength{\parsep}{0pt}
\setlength{\parskip}{0pt}
\begin{enumerate}
  \setlength{\itemsep}{1pt}
  \setlength{\parsep}{0pt}
  \setlength{\parskip}{0pt}
}{\end{enumerate}}

\newenvironment{packed_itemize}{
\setlength{\parsep}{0pt}
\setlength{\parskip}{0pt}
\begin{itemize}
  \setlength{\itemsep}{1pt}
  \setlength{\parsep}{0pt}
  \setlength{\parskip}{0pt}
}{\end{itemize}}

\begin{document} 

\begin{center}
\vspace*{-20pt} 
\centerline{\smalltt INTEGERS: \smallrm ELECTRONIC JOURNAL OF COMBINATORIAL NUMBER THEORY
\smalltt 7 (2007), \#G07} 
\vskip 20pt

\uppercase{\bf The minimum size required of a solitaire army}
\vskip 20pt
{\bf George I. Bell\footnote{\tt 
http://www.geocities.com/gibell.geo/pegsolitaire/}}\\
{\smallit Tech-X Corporation, Boulder, CO 80303, USA}\\
{\tt gibell@comcast.net}\\
\vskip 10pt
{\bf Daniel S. Hirschberg}\\
{\smallit Department of Computer Science, University of California Irvine, Irvine, CA 92697, USA}\\
{\tt dan@ics.uci.edu}\\
\vskip 10pt
{\bf Pablo Guerrero-Garc{\'i}a}\\
{\smallit Dpto. de Matem{\'a}tica Aplicada, Universidad de M{\'a}laga, 29071-M{\'a}laga, Spain}\\
{\tt pablito@ctima.uma.es}\\
\end{center}
\vskip 30pt
\centerline{\smallit Received: 12/20/06, Revised: 5/3/07, Accepted: 7/7/07, Published: 8/10/07}
\vskip 30pt 

\centerline{\bf Abstract}

\noindent
The solitaire army is a one-person peg jumping game
where a player attempts to advance an ``army" of pegs
as far as possible into empty territory.
The game was introduced by John Conway and is also known
as ``Conway's Soldiers".
We consider various generalizations of this game in
different 2D geometries, unify them under a common mathematical framework,
and find the minimum size army capable of advancing a given number of steps.

\pagestyle{myheadings}
\markright{\smalltt INTEGERS: \smallrm ELECTRONIC JOURNAL OF COMBINATORIAL NUMBER THEORY \smalltt 7 (2007), \#G07\hfill}

\thispagestyle{empty} 
\baselineskip=15pt 
\vskip 30pt

\section*{\normalsize 1. Introduction}

\begin{figure}[htbp]
\centering
\epsfig{file=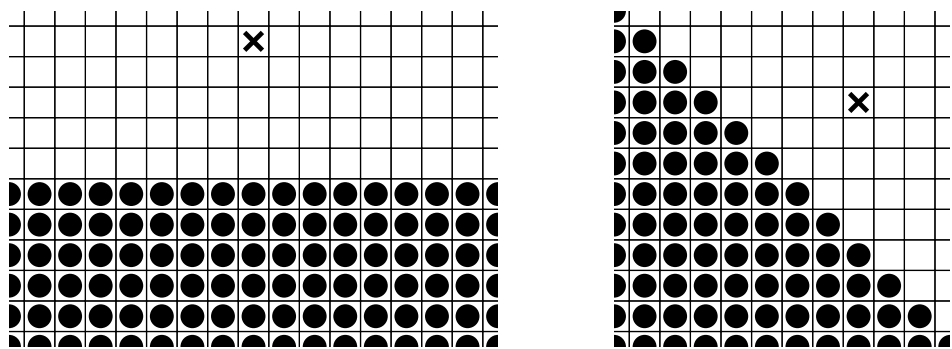}
\caption{(a) Conway's army, and (b) Skew army. X marks the target cell.}
\label{fig1}
\end{figure}

The solitaire army is a one-person game played on an infinite
board \cite{WinningWays}.
Pegs (called men) are placed in all cells below
a certain horizontal line, as in Figure~\ref{fig1}a.
The player can then jump any man over another into an empty cell,
and the jumped man is removed from the board.
The goal is to advance one man as far upward as possible. 
In the usual version of the game, jumps are allowed
only along columns and rows, i.e., \textbf{orthogonal} jumps.

This problem was introduced by John Conway in 1961\footnote{John Beasley personal communication.
The first published appearance of the solitaire army was apparently 1976 \cite{Honsberger}.},
and we refer to this original version as \textit{Conway's army}.
Conway also discovered an elegant proof that no army
restricted to orthogonal jumps can advance more than four steps,
or levels.
His argument using a cleverly chosen weighting function
of board locations has been reproduced many times,
and the puzzle has even found its
way into a best-selling novel \cite{Haddon}.
In addition, the minimum size army to advance 1, 2, 3, or 4
levels was found to be 2, 4, 8, and 20 men, respectively \cite{Beasley}.

The solitaire army problem has been generalized
to include diagonal jumps,
or to occur on a triangular or hexagonal lattice.
This paper has two goals.  The first goal is to summarize
all of these generalizations and unite them under a common framework.
We will discuss only 2D generalizations of the problem,
although it has been studied in higher dimensions as well
\cite{Eriksen95,Eriksen}.
The second goal is to calculate minimum army sizes that are
able to advance a specific number of levels.

\vskip 30pt
\section*{\normalsize 2. Variations on a theme}

One generalization is to take the army in Figure~\ref{fig1}a,
and allow diagonal jumps in either direction
(in addition to orthogonal jumps), giving eight total jump directions possible.
If one plays with this \textit{diagonal army}, it is quite easy
to advance more than four levels.
However, the number of levels the army can advance is still limited,
as we prove in the next section.

Suppose we take the Figure~\ref{fig1}a army, but allow
\textbf{only} diagonal jumps (i.e., orthogonal jumps are not allowed).
Since there are still four jump directions,
one might think that this is equivalent to Conway's army,
but there are subtle differences.
If the jumps can only be diagonal, then they are like
jumps in the game of Checkers.
If we color the board in alternating black and white,
like a Checkers board, we see that the game decomposes into two
separate games that cannot affect one another.
If we consider the game on the color of the finishing cell,
and rotate the board $45^\circ$ (clockwise),
we see that it is related to Conway's army,
but starting from the diagonal or skew front in Figure~\ref{fig1}b.
The army of Figure~\ref{fig1}a with only diagonal jumps
will therefore be called the \textit{skew army}.

\begin{figure}[htbp]
\centering
\epsfig{file=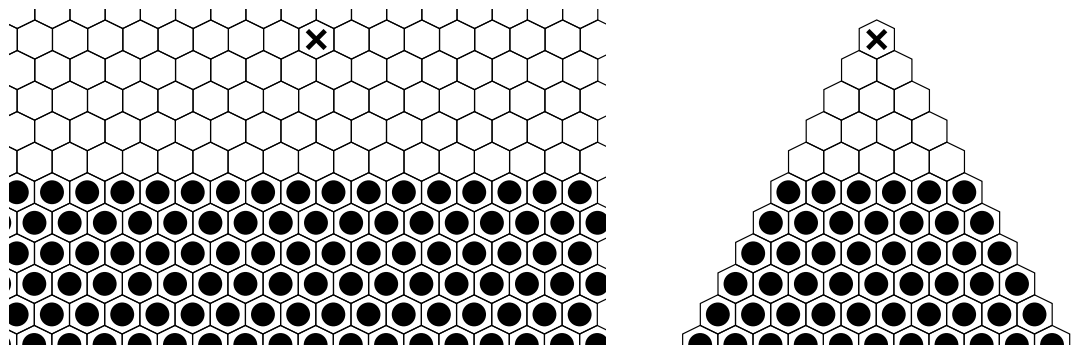}
\caption{(a) Hexagonal army, and (b) Pablito's army (both targeting level 5).}
\label{fig2}
\end{figure}

We could also consider the same problem on the hexagonal lattice in
Figure~\ref{fig2}a, the \textit{hexagonal army}.
Again, our goal is to advance a man as many rows upward as possible,
and we now have six jump directions possible.
Note that playing on a hexagonal board is equivalent to playing
on the usual square lattice board of Figure~\ref{fig1}a,
but allowing only orthogonal jumps plus jumps along one diagonal.

Rather than playing on an infinite hexagonal lattice,
one can play on a triangular board,
with the goal being to finish in the top-most cell (Figure~\ref{fig2}b).
We refer to this version as \textit{Pablito's army}.
It is not clear that this version differs greatly
from the hexagonal army in Figure~\ref{fig2}a,
but we will see that the supply of men is more limited.

In the hexagonal army,
jumps can be made along two diagonal lines
or along a horizontal line.
If a hexagonal army reaches a certain level using \textbf{no}
horizontal jumps at all, then this solution can be
immediately translated into a skew army solution
reaching the same level.
Thus we see that the skew army is very similar to the
hexagonal army, and we will see that solutions can often be
translated in this fashion.

In total, we shall consider five types of solitaire army:
\begin{packed_enumerate}
\item Conway's army, with 4 jump directions.
\item Skew army, with 4 jump directions.
\item Diagonal army, with 8 jump directions.
\item Hexagonal army, with 6 jump directions.
\item Pablito's army, with 6 jump directions.
\end{packed_enumerate}

The diagonal army was first studied by Aigner \cite{Aigner} in 1997,
and has been generalized to $n$-dimensions by Eriksen et al. \cite{Eriksen}.
The hexagonal army was studied by Duncan and Hayes \cite{Duncan} in 1991.
Pablito's army first appeared in 1998 \cite{Pablo1,Pablo2},
and became more widely known after inclusion in a weekly email
puzzle list \cite{Macalester}.
The skew army of Cs{\'a}k{\'a}ny and Juh{\'a}sz \cite{Army} was introduced in 2000.
A 2006 summary article \cite{Niculescu} obtains an upper
bound on the highest level reachable by most of these army types.

\vskip 30pt
\section*{\normalsize 3. Pagoda functions}

As one begins to jump men upward, gaps form in the army.
In order to continue advancing, these gaps must somehow
be spanned, and crossing the gaps introduces more gaps.
It is not hard to believe that the army cannot advance indefinitely,
but proving this is non-trivial.

The idea used by Conway is
to apply a potential or weighting function,
which he called a \textbf{pagoda function}
due to its shape (as a bar graph over the cells of the board).
The \textbf{total weight} of a board position is defined to be
the sum of the weights associated with cells that are occupied by men.
However, not just any weighting function will be useful.
The weighting is called a pagoda function if it has
the property that the total weight \textit{cannot increase}
as men are jumped \cite{WinningWays}.

Suppose we choose $0<\sigma<1$ such that $\sigma^2+\sigma=1$,
i.e., $\sigma=(\sqrt{5}-1)/2\approx 0.618$.
Then
\begin{equation}
\sigma^i+\sigma^{i-1}=\sigma^{i-2}.
\label{SigmaBasic}
\end{equation}

\begin{figure}[htbp]
\centering
\epsfig{file=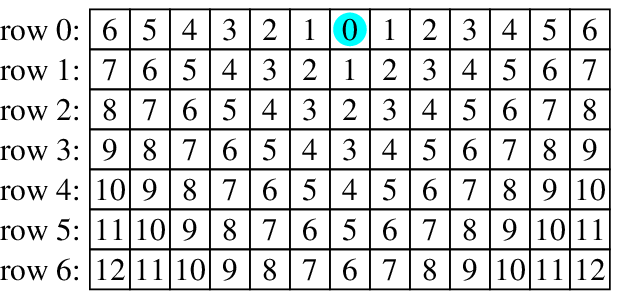}
\caption{Pagoda function exponents for Conway's army.
For a board location labeled $i$ the weighting is $\sigma^i$.}
\label{fig3}
\end{figure}

We assign to each board cell the weight $\sigma^i$,
where the exponent $i$ is the number shown in Figure~\ref{fig3}.
Note that $\sigma^{i+1}<\sigma^i$ since $\sigma<1$.
To set up the solitaire army problem for level $n$,
we place men in all cells of all rows numbered $n$ and higher.
The cell marked ``0" (in light blue) is the target cell we wish to reach,
and has weight $\sigma^0=1$.
Notice that any vertical or horizontal jump into this cell
must use two men with exponents 2 and 1, ending at 0.
Before this jump, the total weight of these three cells was
$\sigma^2+\sigma$, while after the jump it is $\sigma^0=1$
and, by construction, these two quantities are exactly equal.
This weighting is useful because
every orthogonal jump which brings a peg closer to the target cell
maintains the total weight (while any orthogonal
jump away from the target cell decreases the total weight).

If the board location labeled $a$ has weight $\sigma^a$, we
require that our weighting satisfy the \textbf{pagoda condition},
\begin{equation}
\sigma^a+\sigma^b \ge \sigma^c,
\label{PC1}
\end{equation}
for every possible solitaire jump from $a$ over $b$ into $c$.
One possible jump uses the same three cells in the opposite direction,
so we must also have
\begin{equation}
\sigma^c+\sigma^b \ge \sigma^a.
\label{PC2}
\end{equation}
These two conditions can only be met by certain choices of
the exponents $(a,b,c)$.
An exhaustive calculation finds that
three integers $(a,b,c)$ satisfy
Equations~(\ref{PC1}) and (\ref{PC2}) if and only if they
fall into one or more of the following \textbf{types}:
\begin{packed_enumerate}
\item $(a,b,c)$ are drawn from a set of one or two consecutive integers,
or $(a,b,c)$ is any permutation of three consecutive integers.
For example $(3,3,3)$, $(2,2,3)$, $(4,3,4)$,
 $(4,5,3)$ or $(5,4,3)$ all satisfy this criterion.
\item $a=c$ and $a<b$.  For example $(2,5,2)$.
\item $a\ge b$ and $b\le c$. For example $(5,1,3)$.
\end{packed_enumerate}
The first of these types is the most important, and covers
all of the jumps in Figure~\ref{fig3}.
If every jump is (at least) one of these three types, then
the pagoda condition (\ref{PC1}) guarantees that
the total weight \textit{cannot increase} as the game is played.
It follows that if we are to reach the
final board position with weight $\sigma^0=1$,
\textit{the total weight of the starting position must be 
greater than or equal to 1}.

At this point we need a few summation identities in order to calculate
the total weight of our solitaire armies:
\begin{eqnarray}
\sum_{i=n}^{\infty} \sigma^i & = & \frac{\sigma^n}{1-\sigma} = \sigma^{n-2}, \label{eq:sum1} \\
\sum_{i=n}^{\infty} \sigma^{2i} & = & \frac{\sigma^{2n}}{1-\sigma^2} = \sigma^{2n-1}, \label{eq:sum1a} \\
\sum_{i=n}^{\infty} i\sigma^i & = &
n\sum_{i=n}^{\infty} \sigma^i+\sum_{i=n+1}^{\infty}\sum_{j=i}^{\infty} \sigma^j=
n\sigma^{n-2} +\sigma^{n-3}. \label{eq:sum2}
\end{eqnarray}

Any Conway's army capable of reaching the final cell must have total
weight greater than or equal to $1$.
We can compute the total weight of row $n$, assuming that it
is entirely filled by men, using the Identities (\ref{eq:sum1}) and (\ref{SigmaBasic}), as
$$R_{1,n} = \sum_{i=n}^{\infty} \sigma^i + \sum_{i=n+1}^{\infty} \sigma^i
= \sigma^{n-2}+\sigma^{n-1} = \sigma^{n-3},$$
where the first subscript refers to the army type, and the second the row.
Therefore, the sum of rows $n$ and beyond is
$$S_{1,n} = \sum_{i=n}^{\infty} R_{1,i} = \sum_{i=n}^{\infty} \sigma^{i-3} = \sigma^{n-5}.$$
The total weight of rows $5$ and higher is $S_{1,5}=\sigma^0=1$,
so it is impossible for any \textit{finite} army to reach level 5.
This completes the proof that a Conway's army cannot
reach a level greater than 4.

\begin{figure}[htbp]
\centering
\epsfig{file=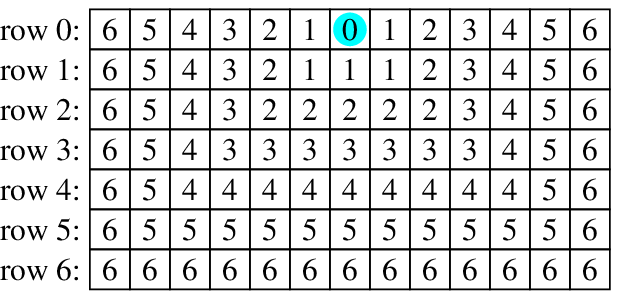}
\caption{Pagoda function exponents for skew and diagonal armies.}
\label{fig4}
\end{figure}

The exponents of Figure~\ref{fig3} are not a valid pagoda function
when diagonal jumps are allowed,
this requires a different weighting function
with exponents as shown in Figure~\ref{fig4}.
Note that both of these exponent patterns are based on a
distance metric from the target cell,
where Figure~\ref{fig3} uses the Manhattan or taxi cab metric
($\ell_1$ norm),
and Figure~\ref{fig4} the number of chess king moves
($\ell_\infty$ norm).

For the skew army, recall that, for each row,
only every other cell contains a man.
Thus the sum of row $n$, using (\ref{eq:sum1a}) is
$$R_{2,n} = (n+1)\sigma^n + 2\sigma^n\sum_{i=1}^{\infty} \sigma^{2i} =
(n+1)\sigma^n + 2\sigma^{n+1}.$$
The sum of rows $n$ and beyond is, using~(\ref{eq:sum1}) and (\ref{eq:sum2}),
$$S_{2,n} = \sum_{i=n}^{\infty} R_{2,i} = \sigma^{n-3}((n-1)\sigma+3).$$
Direct calculation shows that $S_{2,6}\approx 1.44$,
but that $S_{2,7}<1$, so level 7 cannot be reached.
At this point we do not know whether level 6 can be reached,
but we will show that it is possible.

For the diagonal army, the sum of row $n$
using (\ref{eq:sum1}) and (\ref{SigmaBasic})  is
$$R_{3,n} = (2n+1)\sigma^n + 2\sum_{i=n+1}^{\infty} \sigma^i =
2n\sigma^n + \sigma^{n-3}.$$
Hence the sum of rows $n$ and beyond is, using (\ref{eq:sum1}) and (\ref{eq:sum2}),
$$S_{3,n} = \sum_{i=n}^{\infty} R_{3,i} = \sigma^{n-5}((4n-2)\sigma+3-2n).$$
Direct calculation shows that $S_{3,8}\approx 1.31$,
but that $S_{3,9}<1$, so level 9 cannot be reached.
Again, we will show that level 8 can be reached.

\begin{figure}[htbp]
\centering
\epsfig{file=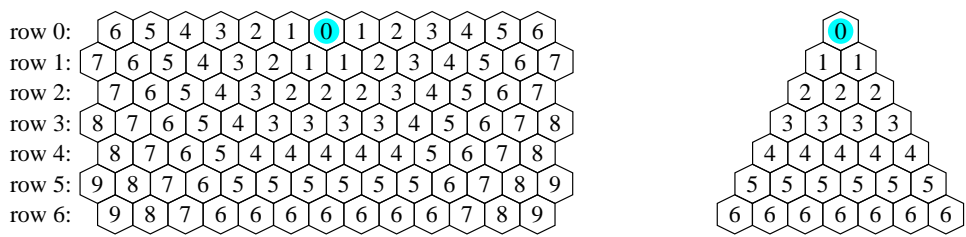}
\caption{Pagoda function exponents for (a) hexagonal army, and (b) Pablito's army.}
\label{fig5}
\end{figure}

For Pablito's army, we use the weighting function given by
Figure~\ref{fig5}b.  The row sum is simply
$$R_{5,n} = (n+1)\sigma^n,$$
and so the sum of rows $n$ and beyond is
$$S_{5,n} = \sigma^{n-3}((n+1)\sigma+1).$$
Here $S_{5,6}\approx 1.26$, but $S_{5,7}<1$,
and 6 is the highest level that can be reached.
The hexagonal army uses the weighting function of Figure~\ref{fig5}a,
which adds on a bit more
$$R_{4,n} = R_{5,n} + 2\sum_{i=n+1}^{\infty} \sigma^i=(n+1)\sigma^n+2\sigma^{n-1},$$
and, using~(\ref{eq:sum1}) again,
$$S_{4,n} = S_{5,n}+2\sigma^{n-3}=\sigma^{n-3}((n+1)\sigma+3).$$
One can calculate that $S_{4,7}\approx1.16$, but $S_{4,8}<1$, so the hexagonal
army cannot reach level 8.
Table~\ref{table1} summarizes upper bounds on the number of
levels each type of army can advance.
We will show that all these upper bounds can be attained.

\begin{table}[htbp]
\begin{center} 
\begin{tabular}{ | c | c | c | c | c | c | }
\hline		
Conway's & Skew & Pablito's & Hexagonal & Diagonal  \\
\hline
\hline
4 & 6 & 6 & 7 & 8 \\
\hline
\end{tabular}
\caption{Upper bounds on the highest level reachable by each type of army.} 
\label{table1}
\end{center} 
\end{table}

\vskip 30pt
\section*{\normalsize 4. Minimum size armies}

\subsection*{\normalsize 4.1 A simple lower bound}

The pagoda function arguments give an upper bound
on the number of levels an army can reach,
and they can also give an idea of what an army
capable of reaching a certain level must look like.
Instead of adding the pagoda function by rows, we can add it in
order of decreasing weight (increasing exponent) and see exactly
when this sum becomes at least 1.
This provides a lower bound on the army size capable
of reaching a certain level.

For example, for Conway's army to reach level 4,
the front-most men are in row 4 in Figure~\ref{fig3},
and one has weight $\sigma^4$.
After that there are 3 with weight $\sigma^5$, 5 with weight $\sigma^6$, etc.
The total weight of the army in order of decreasing weight is
$$\sigma^4+3\sigma^5+5\sigma^6+7\sigma^7+\cdots$$
When does this sum become at least $1$?
If we start with 19 men, the greatest possible total weight is
$$T=\sigma^4(1+3\sigma+5\sigma^2+7\sigma^3+3\sigma^4).$$
We can calculate this by collapsing the sum from right to left using the
Identity (\ref{SigmaBasic}).
This can be accomplished in tableau form,
which is really a form of peg solitaire in 1-dimension:

\begin{tabular}{ r r r r r }
$1$ & $3$ & $5$ & $7$ & $3$ \\
 & & $3$ & $-3$ & $-3$ \\
 & $4$ & $-4$ & $-4$ & \\
$4$ & $-4$ & $-4$ & & \\
\hline
$5$ & $3$ & $0$ & $0$ & $0$ \\
\end{tabular}

Therefore $T=\sigma^4(5+3\sigma)$.
In fact $T$ is identically equal to 1 because $5+3\sigma=\sigma^{-4}$.
We can show this by deriving an identity which
connects the golden ratio $\sigma^{-1}$ and Fibonacci numbers:
\begin{equation}
\sigma^i=(-1)^i[F_{i-1}-F_i\sigma], \hspace{.5in} i\in\mathbb{Z},
\label{Fibonacci}
\end{equation}
where $F_i$ are the Fibonacci numbers, identified by $F_1=F_2=1$, and
$F_i=F_{i-2}+F_{i-1}$.
Equation~(\ref{Fibonacci}) can be proved by induction, and applies
to all $i\in\mathbb{Z}$ if we extend the Fibonacci numbers by defining
$F_0=0$, $F_{-i}=(-1)^{i+1}F_i$.
Table~\ref{table2} gives powers of $\sigma$
in terms of linear combinations of $1$ and $\sigma$.

\begin{table}[htbp]
\begin{center} 
\begin{tabular}{ | c | c | c | }
\hline
$n$ & $\sigma^n$ & $\sigma^{-n}$ \\		
\hline
\hline
1 & $\sigma$ & $1+\sigma$ \\
\hline
2 & $1-\sigma$ & $2+\sigma$ \\
\hline
3 & $-1+2\sigma$ & $3+2\sigma$ \\
\hline
4 & $2-3\sigma$ & $5+3\sigma$ \\
\hline
5 & $-3+5\sigma$ & $8+5\sigma$ \\
\hline
6 & $5-8\sigma$ & $13+8\sigma$ \\
\hline
7 & $-8+13\sigma$ & $21+13\sigma$ \\
\hline
8 & $13-21\sigma$ & $34+21\sigma$ \\
\hline
9 & $-21+34\sigma$ & $55+34\sigma$ \\
\hline
10 & $34-55\sigma$ & $89+55\sigma$ \\
\hline
\end{tabular}
\caption{Powers of $\sigma=(\sqrt{5}-1)/2$.} 
\label{table2}
\end{center} 
\end{table}

This shows that a lower bound on the army size to reach level 4 is 19.
Using the same argument on pagoda functions for each type of army,
we can derive lower bounds on army size for
all five armies over all feasible levels.
However, it turns out that the above pagoda functions tend to overestimate
how well we can utilize the men towards the
left or right edges of the army.
An improved lower bound can be found using a modification
of the pagoda functions above, as we show in the next section.

If $L_n$ is the minimum size army needed to reach level $n$ in some geometry,
the final jump must involve two pegs at levels $n-1$ and $n-2$.
Therefore, it is clear that
\begin{equation}
L_n \ge L_{n-1}+L_{n-2} ,
\label{Fibonacci1}
\end{equation}
showing that if we have lower bounds on $L_{n-1}$ and $L_{n-2}$,
we can obtain lower bounds for all higher $L_n$.
Moreover, since $L_0=1=F_2$ and $L_1=2=F_3$,
$L_n$ in any geometry is always bounded below
by the Fibonacci number $F_{n+2}$,
\begin{equation}
L_n \ge F_{n+2} .
\label{GlobalLowerBound}
\end{equation}

\subsection*{\normalsize 4.2 An improved lower bound}

\begin{figure}[htbp]
\centering
\epsfig{file=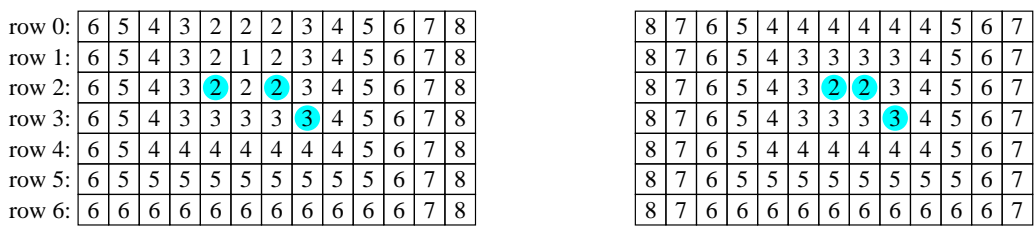}
\caption{Improved pagoda function exponents for the diagonal or skew armies
when (a) the last two jumps are diagonal, and (b) the last jump is vertical.}
\label{fig6}
\end{figure}

Let us consider the diagonal army and suppose the
last two jumps are diagonal jumps.
Before these jumps are made, the board must look
like Figure~\ref{fig6}a (or its reflection),
with three men shown in light blue.
This figure shows exponents for an improved pagoda function,
with the total weight of the position shown
$2\sigma^2+\sigma^3=\sigma+\sigma^2=1$.
However, the lower bound on the army size
will be tighter using this pagoda function.
This is because the weight of each row
in Figure~\ref{fig6}a is less than the weight of that row
in Figure~\ref{fig4}---for example row 5 contains nine 5's,
while in Figure~\ref{fig4} it had eleven 5's.
The alert reader may notice that the value of the final
position is now $\sigma^2$, but this is not significant
because we still need to pass through a board position with
total weight $1$.

As an example let us compute the lower bound for the diagonal army
to reach level 5 using the pagoda function in Figure~\ref{fig6}a.
We can place at most 9 men in row 5 with exponent 5, and
13 men in rows 5 and 6 with exponent 6.
We can calculate (using Table~\ref{table2}) that
$9\sigma^5+3\sigma^6=1-\sigma^8<1$,
but $9\sigma^5+4\sigma^6=1+\sigma^7>1$, which
gives a lower bound on the army size of $9+4=13$.

What if the last two jumps are not diagonal?
If the last jump is vertical, then we must have the situation of
Figure~\ref{fig6}b (or its reflection).
This configuration will give a larger lower bound,
because there are only eight 5's in row 5.
The same pagoda function works if the last jump is diagonal
and the second to the last jump is vertical.
This argument proves that at least 143 men are required to
reach level 8 if either of the last two jumps is vertical.
But using two diagonal jumps as in Figure~\ref{fig6}a,
we obtain a lower bound of $15+19+23+27+31+7=122$.
Assuming these lower bounds are close to the actual minimums,
the smallest army capable of reaching level 8
must finish with two diagonal jumps.

\begin{figure}[htbp]
\centering
\epsfig{file=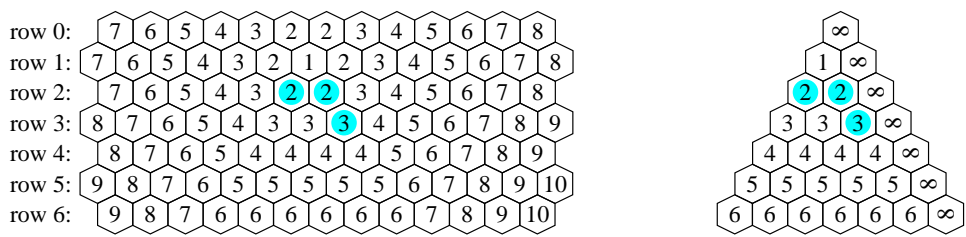}
\caption{Improved pagoda function exponents for (a) hexagonal army, and (b) Pablito's army.}
\label{fig7}
\end{figure}

Similar reasoning on the hexagonal army gives the
improved pagoda function exponents in Figure~\ref{fig7}a.
For Pablito's army, the best pagoda function is shown in Figure~\ref{fig7}b.
The exponent values $\infty$ are just board weights of $\sigma^\infty=0$.
This pagoda function indicates that one edge of the board
is of no use towards accomplishing our goal;
we might as well leave it empty.

If we take each of these improved pagoda functions, and calculate the
total weight in order of increasing exponents, we obtain the lower bounds
on army size given in Table~\ref{table3}.
The reader is encouraged to construct maximum total weight armies
as in the level~5 example to understand the bounds in
Table~\ref{table3}.

\begin{table}[htbp]
\begin{center} 
\begin{tabular}{ | c | c | c | c | c | c | c | }
\hline
Level \# & \multicolumn{6}{|c|}{Lower bounds on the army size to reach level $n$}\\		
($n$) & $F_{n+2}$ & Conway's & Skew & Pablito's & Hexagonal & Diagonal\\
\hline
\hline
1 & 2 & 2 & 2 & 2 & 2 & 2\\
\hline
2 & 3 & 4 & 3 & 3 & 3 & 3\\
\hline
3 & 5 & 8 & 5 & 5 & 5 & 5\\
\hline
4 & 8 & 19 & 9 & 9 & 9 & 8\\
\hline
5 & 13 & Impossible & 18 & 18 & 16 & 13\\
\hline
6 & 21 & Impossible & 43 & 51 & 35 & 23\\
\hline
7 & 34 & \multicolumn{3}{|c|}{Impossible} & 140 & 45\\
\hline
8 & 55 & \multicolumn{4}{|c|}{Impossible} & 122\\
\hline
9 & 89 & \multicolumn{5}{|c|}{All Impossible}\\
\hline
\end{tabular}
\caption{Lower bounds on the army size required to reach level $n$,
from equation~(\ref{GlobalLowerBound}) and
the pagoda function of Figures~\ref{fig3}, \ref{fig6}a, \ref{fig7}b,
\ref{fig7}a and \ref{fig6}a respectively.}
\label{table3}
\end{center} 
\end{table}

\subsection*{\normalsize 4.3 Minimum size diagonal armies}

How good are the lower bounds in Table~\ref{table3}?
It turns out they are very good---to see this
we must produce small size diagonal armies capable
of reaching a certain level.
One way to do this is in an inductive fashion.
It is easy to find an 8-man diagonal army capable of reaching level 4.
Next, we try to find an army at level 5 that, after some sequence
of jumps, reproduces the previous 8-man army one row forward from
the starting line.
This technique has been used \cite{Honsberger,Aigner} for
Conway's armies, but in general it is useful only
for the lowest levels.
To find the smallest army that can reach levels 7 or 8,
a different technique is needed.

\begin{figure}[htb]
\centering
\epsfig{file=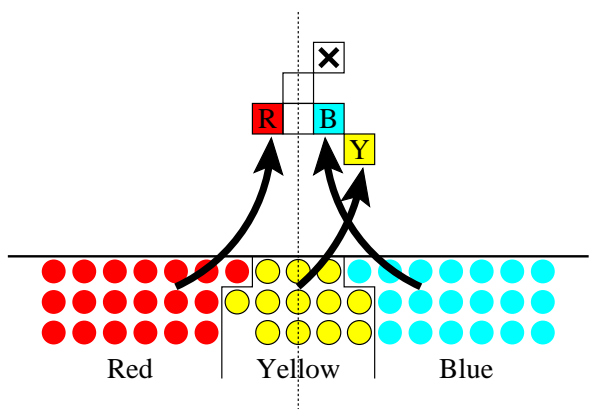}
\caption{An army divided into 3 regiments.}
\label{fig8}
\end{figure}

A better idea comes from our target configuration of three men in
Figure~\ref{fig6}a, reproduced in Figure~\ref{fig8}.
Following Duncan and Hayes \cite{Duncan}, we separate our army into
differently colored \textit{regiments}, and the goal of each regiment
is to reach the target cell of its color.
We may even be able to use symmetric regiments (red and blue),
symmetric about the dashed line,
with the remaining yellow target cell filled by the
(nonsymmetric) central yellow regiment.
We cannot always achieve such a symmetric army,
but we shall see that often it is possible.

\begin{table}[htbp]
\begin{center} 
\begin{tabular}{ | c | l | c | c | }
\hline
Type & Exponents of the three jump cells & Jump right loses & Jump left loses \\		
\hline
\hline
1.1.1 & $(a,a,a)$ & $\sigma^a$ & $\sigma^a$ \\
\hline
1.2.1 & $(a,a,a-1)$ & $\sigma^{a+2}$ & $\sigma^{a-1}$ \\
\hline
1.2.2 & $(a,a-1,a)$ & $\sigma^{a-1}$ & $\sigma^{a-1}$ \\
\hline
1.2.3 & $(a,a-1,a-1)$ & $\sigma^a$ & $\sigma^{a-1}+\sigma^{a+1}$ \\
\hline
1.2.4 & $(a-1,a,a-1)$ & $\sigma^a$ & $\sigma^a$ \\
\hline
1.3.1 & $(a,a-1,a-2)$ & $0$ & $2\sigma^{a-1}$ \\
\hline
1.3.2 & $(a,a-2,a-1)$ & $2\sigma^a$ & $2\sigma^{a-1}$ \\
\hline
1.3.3 & $(a-1,a,a-2)$ & $0$ & $2\sigma^a$ \\
\hline
\hline
2.1.1 & $(a-i,a,a-i)$, any $i>0$ & $\sigma^a$ & $\sigma^a$ \\
\hline
\hline
3.1.1 & $(a,a-i-j,a-j)$, any $i>0$ and $j>0$ & $\ge\sigma^a$ & $\ge\sigma^{a-j}$ \\
\hline
\end{tabular}
\caption{The amount lost by various types of jump
(these types further subdivide the types presented in section~3).
For each type, $a$ is the largest exponent.} 
\label{table4}
\end{center} 
\end{table}

{\bf Proposition 1} Let us consider an army on some board,
with a pagoda function weighting $\sigma^e$
where all the exponents $e$ are integers.
Suppose there are no jumps of Type 1.2.1 in Table~\ref{table4},
and that the maximum exponent of any occupied cell is $E$.
Then
\begin{packed_enumerate}
\item A jump either loses nothing, or loses \textit{at least} $\sigma^E$.
\item If we perform any number of jumps that lose nothing,
the maximum exponent $E$ of any occupied cell cannot increase.
\end{packed_enumerate}

{\it Proof:} The proof is contained in Table~\ref{table4},
which shows the amount lost by all possible jump exponents
in a valid pagoda function of the form $\sigma^e$.
Each line lists one pattern of three exponents in a jump,
and the maximum exponent is always $a$. 
The amount lost is always either zero or $\ge\sigma^a$
(with the exception of Type 1.2.1, which by hypothesis cannot occur).
Note that none of the pagoda functions in this paper
contain a jump of Type 1.2.1.
For the second part of the proposition,
note that the two jump types that
lose nothing always place men in cells with smaller exponents.

\begin{figure}[htb]
\centering
\epsfig{file=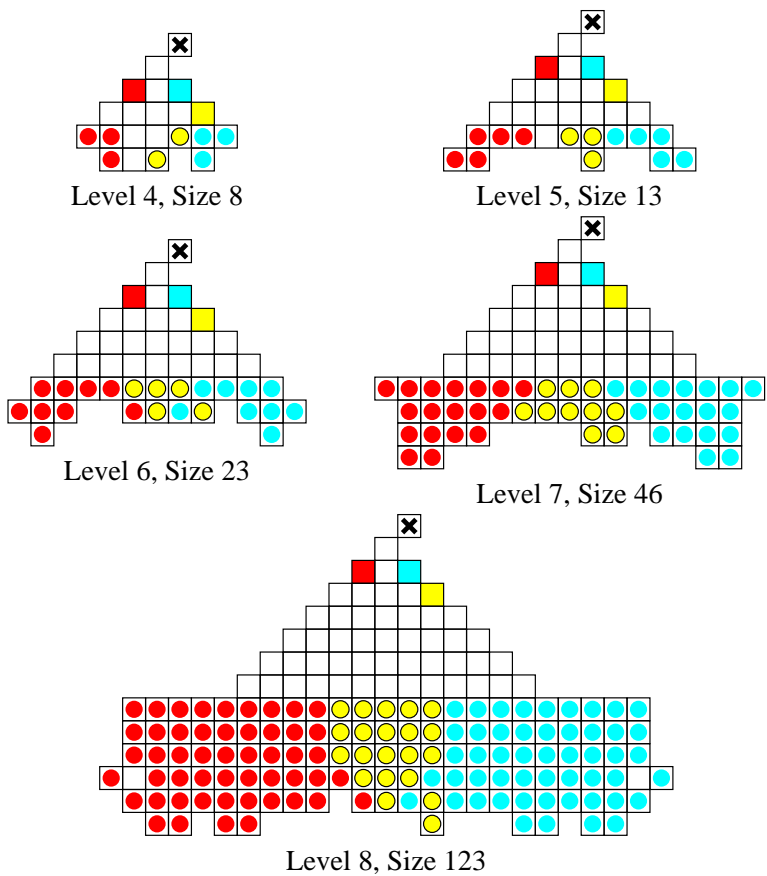}
\caption{Minimum size diagonal armies for levels 4--8.}
\label{fig9}
\end{figure}

{\bf Theorem 1} The armies of Figure~\ref{fig9}
have the smallest possible size for the levels they attain.

{\it Proof:} The level 4--6 armies have size equal to the
lower bound in Table~\ref{table3},
hence they must be of minimum size.
The armies to levels 7 and 8 have size only one larger than the lower bound
of Table~\ref{table3}, so we need to show that no 45-man army can
advance 7 levels, and that no 122-man army can advance 8 levels.
For any configuration of men, the \textbf{slack} is the amount by which the
total weight exceeds 1.
Any army capable of reaching a certain level must begin with
slack greater than or equal to zero.
The \textbf{target board position} is the
3-man position of Figure~\ref{fig6}a,
plus any other men that happen to be left.
If only these three men remain,
the slack of the target board position is identically zero,
otherwise it could be positive.
The basic argument is that if the slack ever becomes negative,
the target board position can no longer be reached.
The army of 45 men with the greatest total weight has slack
\begin{equation}
S=\sigma^7\left(13+17\sigma+15\sigma^2\right)-1
=\sigma^7\left(7-11\sigma\right)=\sigma^{11}+\sigma^{13}
\label{army45}
\end{equation}
using Table~\ref{table2}.

We note that the starting army has maximum exponent $E=9$.
We have from Proposition~1 that any jump either
loses nothing or loses at least $\sigma^9$,
which would make the slack negative because $\sigma^9 > \sigma^{11}+\sigma^{13}$.
Even if we perform any number of jumps that lose nothing,
$E$ cannot increase,
and so any jump that loses weight will make the slack negative.
The only remaining possibility is that we could reach some
target board position using \textit{only} jumps that lose nothing.
Since the starting slack is $\sigma^{11}+\sigma^{13}$, the target
board position must also have this slack, and the only
way this can happen is if we have the 3-man position of Figure~\ref{fig6}a,
plus at least one man at row 11 or higher.
However, the starting army has no men below row 9, so getting a man to
row 11 or higher requires a downward jump, which loses total weight.
This contradicts our assumption that the solution uses only
jumps that lose nothing.

This only proves that a 45-man army with maximum total weight
cannot reach level 7.
It still could be possible that a 45-man army with somewhat smaller
total weight (but still at least 1) could reach level 7.
If we take the starting position with maximum total weight
and move a man from a cell with exponent $n$ to a cell with exponent $n+1$,
the slack decreases by $\sigma^{n+2}$
(this follows directly from (\ref{SigmaBasic}) with $i=n+2$).
So if we move a man from exponent $9$ to $10$ then
$S=\sigma^{13}$.
Proposition~1 again implies that the slack cannot be reduced to zero,
as the highest exponent $E=10$.
If we move any other man to a cell with higher exponent
the slack becomes negative,
so a 45-man army cannot reach level 7.

We use the same argument for a 122-man army to reach level 8,
but there are more cases.
The maximum total weight 122-man army has slack
\begin{equation}
S=\sigma^8\left(15+19\sigma+23\sigma^2+27\sigma^3+31\sigma^4+7\sigma^5\right)-1
=\sigma^{14}+\sigma^{16}.
\label{army122}
\end{equation}
This army has $E=13$, and again the slack cannot be reduced to zero.
If we let $n_i$ be the number of men with exponent $i$, the possible cases
are listed in Table~\ref{table5}.
In cases 1--6, the slack is too small to reduce to zero by Proposition~1.
Case 7 has zero slack but there is no way to get to the man
with exponent 17 without losing total weight, and this man must be
used since the slack of the army without him is negative.

\begin{table}[htbp]
\begin{center} 
\begin{tabular}{ | c | c | c | c | c | c | c | c | }
\hline	
Case & $n_{12}$ & $n_{13}$ & $n_{14}$ & $n_{15}$ & $n_{16}$ & $n_{17}$ & Slack \\
\hline
\hline
1 & 31 & 7 & 0 & 0 & 0 & 0 & $\sigma^{14}+\sigma^{16}$ \\
\hline
2 & 31 & 6 & 1 & 0 & 0 & 0 & $2\sigma^{16}$ \\
\hline
3 & 31 & 6 & 0 & 1 & 0 & 0 & $\sigma^{16}$ \\
\hline
4 & 30 & 8 & 0 & 0 & 0 & 0 & $\sigma^{16}$ \\
\hline
5 & 31 & 5 & 2 & 0 & 0 & 0 & $\sigma^{18}$ \\
\hline
6 & 31 & 6 & 0 & 0 & 1 & 0 & $\sigma^{18}$ \\
\hline
7 & 31 & 6 & 0 & 0 & 0 & 1 & 0 \\
\hline
Figure~\ref{fig9}e & 28 & 11 & 0 & 0 & 0 & 0 & 0 \\
\hline
\end{tabular}
\caption{Number of men of each exponent,
and slack, for an army of 122 men,
or 123 men (last row).
In all cases
$n_i=4i-17$ for $8\le i\le 11$.} 
\label{table5}
\end{center} 
\end{table}

It remains, of course, to show that the armies of Figure~\ref{fig9}
are capable of reaching their specified levels.
For the smallest armies, the jumps leading to a solution
are relatively obvious.
But for the armies to reach levels 7 and 8,
the jumps may not be so clear.
These make interesting puzzles to solve by hand,
but see \cite{Bellweb} to view solutions.

All the armies in Figure~\ref{fig9} have zero slack,
or equivalently total weight exactly 1.
We can see this by calculating the weight of the
entire army, or each regiment.
For example, for the 123-man army to reach level~8,
the total weight of the blue or red
51-man regiments is
$T=\sigma^8(5+7\sigma+9\sigma^2+12\sigma^3+13\sigma^4+5\sigma^5)$.
We can simplify this using the tableau:

\begin{tabular}{ r r r r r r }
$5$ & $7$ & $9$ & $12$ & $13$ & $5$ \\
 & & & $5$ & $-5$ & $-5$ \\
 & & $8$ & $-8$ & $-8$ & \\
 & $9$ & $-9$ & $-9$ & & \\
$8$ & $-8$ & $-8$ & & & \\
\hline
$13$ & $8$ & $0$ & $0$ & $0$ & $0$ \\
\end{tabular}

By Table~\ref{table2}, $13+8\sigma=\sigma^{-6}$,
so $T=\sigma^2$, exactly the weight
of the target cell.
When playing from an army with zero slack,
every jump must be of the Type 1.3.1, $(a,a-1,a-2)$,
in other words over a decreasing sequence of
exponents\footnote{Technically the jump could also
be of Type 1.3.3, namely $(a-1,a,a-2)$, but this pattern never occurs
in the pagoda functions in this paper.}.
For the diagonal armies of Figure~\ref{fig9},
this means that nearly all jumps must be upward,
and that horizontal jumps can only occur near
the right and left edges of the army.

Minimum size armies \textit{usually} have zero slack
(with respect to the improved pagoda functions).
Intuitively, if an army with maximum exponent $E$
does not have zero slack, then Proposition~1
indicates it must have a slack of at least $\sigma^E$,
in which case its size could probably be reduced by 1.
However, this is not always the case, and in Section~5.1
we will see an example of a minimum size army
with non-zero slack.

\subsection*{\normalsize 4.4 Other minimum size armies}

Can we use the same arguments as in the previous section
to find minimum size armies of other types?
Suppose we look at Pablito's army to reach level 5, here
the maximum total weight 18-man army has slack
\begin{equation}
S=\sigma^5\left(5+6\sigma+7\sigma^2\right)-1=2\sigma^9.
\label{army18}
\end{equation}
Since the maximum exponent $E=7$, and $S=2\sigma^9<\sigma^7$,
this army cannot reduce its slack to zero, by Proposition~1.
However, if we move two starting men from a cell with exponent 7
to a cell with exponent 8, the slack becomes zero.
What this indicates is that an 18-man army capable
of reaching level 5 probably has 5, 6, 5, and 2 men
with exponents 5, 6, 7, and 8, respectively.
If you look for such an army, however, you
will not be able to find one that can reach level 5.
This is also a problem with Conway's army to
level 4---we saw at the beginning of Section~4 that there
are 19-man armies with zero slack,
yet none of them can reach level 4 \cite{Beasley,Aigner}.
For these situations, we need to develop a new technique
to prove our armies have minimum size.

\vskip 30pt
\section*{\normalsize 5. An integer programming model}

\subsection*{\normalsize 5.1 The model}

We can generalize the solitaire army problem by allowing any
integer number of men in a cell.
A jump adds $(-1,-1,+1)$ to a consecutive triple of cells.
This reduces the problem to a question in linear algebra.
We can even state the question of minimum army size directly
as an integer programming (IP) model.
Suppose we have a board
with $M$ cells, and let $N$ be the total number of
jumps possible on this board.

Parameters are:
\begin{packed_itemize}
\item Let $A_{m,n}$ be the jump matrix describing
the effect of each jump on the number of men in each cell.
For example if jump 3 starts with a man at cell 7,
jumps over a man at cell 8 and ends at cell 9,
we would have $A_{7,3}=-1$, $A_{8,3}=-1$ and $A_{9,3}=+1$,
with all other $A_{m,3}=0$ for $m\not\in\{7,8,9\}$.
\item Let $\mbox{SMAX}_m$ be the maximum number of starting men
in cell $m$.
For the solitaire army problem, we would set $\mbox{SMAX}_m$
to $0$ up to a certain row, and $1$ beyond it.
\item Let $\mbox{SMIN}_m$  be the minimum number of starting men
in cell $m$.
Normally we let this be zero everywhere unless we want to
force a man at some location.
\item Let $\mbox{FIN}_m$ be the final number of men in cell $m$.
For a solitaire army problem,
we would set all these to zero except for the final top-most cell,
set to 1.
\end{packed_itemize}

Decision variables are:
\begin{packed_itemize}
\item $\mbox{STA}_m$, equal to 1 if there is initially a man in cell $m$,
0 otherwise.  A binary variable.
\item $J_n$, the number of jumps of index $n$ present in the solution.
A non-negative integer variable.
\end{packed_itemize}

Objective function is to:

$$\mbox{minimize total starting men} = \sum_{m\in M} \mbox{STA}_m.$$

Constraints that must be satisfied by the solution:

\begin{tabular}{ l l }
$\mbox{STA}_m+\sum_{n\in N}J_n A_{m,n} = \mbox{FIN}_m, \forall m\in M$, &
jumps applied to the starting configura-\\
& tion gives the final configuration, and \\
 & \\
$\mbox{SMIN}_m \le \mbox{STA}_m \le \mbox{SMAX}_m, \forall m\in M$, &
starting number of men must fall\\
& within bounds.
\end{tabular}

It is important to realize that the above IP
problem is not equivalent to the original solitaire army
problem\footnote{It is possible to create an IP model which is equivalent
to the original peg solitaire problem \cite{Kiyomi,IPModel2},
but this model will be considerably more complicated,
and take much longer to solve optimally.}.
In the original problem the jumps occur in a specific order,
and it may be impossible to order the jumps in a solution to the
IP model so that any cell contains either zero or one man at any
time during the solution, and each jump can be legally executed.
See \cite{Bell} for examples of problems where the IP model is solvable
but the solution cannot be translated into a peg solitaire solution.

However, any solution to the original solitaire army
problem is a solution to the IP model.
Therefore, the number coming out of the above IP minimization
is a lower bound on the size of a solitaire army.
In practice we have never seen a case where we cannot order the jumps
from a solution to the IP model to create a solution to the
solitaire army problem, although this ordering process can be quite
difficult for the largest problems.

The IP model can only be satisfied by a linear
combination of jumps that takes the board from the initial state
to the final state.
This is related to the \textit{lattice criterion}
that appears in the literature \cite{DezaComb,DezaLattice}.
The difference is that the lattice criterion
allows $J_n$ to take on any integer value, while
the above IP model requires $J_n\ge 0$.
A jump $J_n=-1$ changes three consecutive cells by $(+1,+1,-1)$,
which \textit{adds} a man rather than removing one.
If we allow $J_n$ to take on any integer value,
a one-man army can advance an arbitrary number of levels.
This is easy to see using the tableau:

\begin{tabular}{ r r r r }
$+1$ & $0$ & $0$ & $0$\\
$-1$ & $+1$ & $+1$ & \\
 & $-1$ & $-1$ & $+1$\\
\hline
$0$ & $0$ & $0$ & $+1$\\
\end{tabular}

\begin{table}[htbp]
\begin{center} 
\begin{tabular}{ | c | c | c | c | c | }
\hline
 & \multicolumn{4}{|c|}{Lower bounds on the army size to reach level $n$} \\		
Level \# ($n$) & Conway's & Skew & Pablito's & Hexagonal \\
\hline
\hline
4 & 20 (1 sec) & $-$ & $-$ & $-$ \\
\hline
5 & Impossible & 19 (1 sec) & 19 (1 sec) & 17 (1 sec) \\
\hline
6 & Impossible & 46 (2 sec) & 53 (7 sec) & 36 (3 sec) \\
\hline
7 & \multicolumn{3}{|c|}{Impossible} & 144 (Section~5.2) \\
\hline
\end{tabular}
\caption{Lower bounds on the army size required to reach level $n$,
using the IP model.  Solve times in parenthesis.}
\label{table6}
\end{center} 
\end{table}

The only difficulty in applying the IP model is in deciding
how large to make the board.
If we make the board too small, a solution may be eliminated,
giving us a larger army than the smallest possible one.
The pagoda function arguments can be useful here to decide 
how many rows to include.

\begin{figure}[htb]
\centering
\epsfig{file=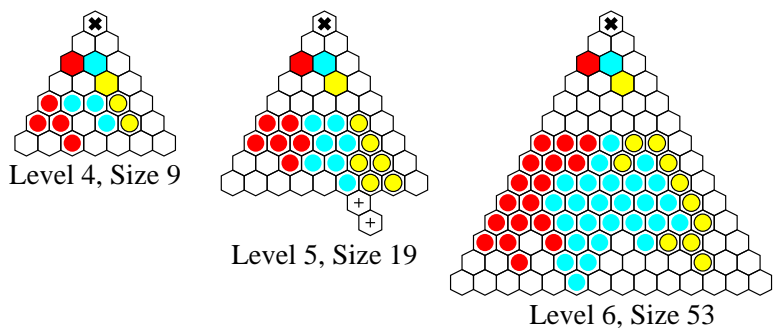}
\caption{Minimum size Pablito's armies to reach levels 4--6.}
\label{fig10}
\end{figure}

We have programmed the IP models in the GAMS language,
and solved them using several commercial solvers available at
the NEOS server \cite{NEOS}.
Table~\ref{table6} shows the results of the IP model applied to
all remaining unsolved problems.
All bounds shown have been improved from those given in Table~\ref{table3}.
Xpress-MP solved all the problems in under 10 seconds,
with the exception of the hexagonal army to level 7,
to be discussed in the Section~5.2.
All solve times quoted are for the Xpress-MP
solver run on the NEOS server {\tt newton.mcs.anl.gov},
a 2.4GHz Intel Xenon with 2 GB of memory.

\begin{figure}[htb]
\centering
\epsfig{file=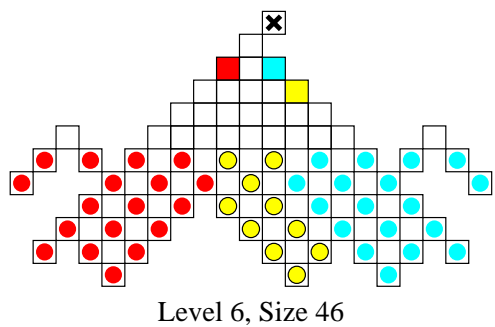}
\caption{A minimum size skew army to reach level 6.}
\label{fig12}
\end{figure}

For all cases except for the hexagonal army to level 7,
we can find armies with size equal to
the lower bounds in Table~\ref{table6}.
Figure~\ref{fig10} shows minimum size Pablito's armies.
The Pablito's army to reach level 5 is the only army shown in
this paper with non-zero slack.
We can convert this army into another minimum size army
with zero slack by moving the bottom two yellow men to the
cells marked by crosses.
These zero slack armies only use upward jumps
in their solution.
Hence these solutions can be immediately converted to
skew army solutions, and those for levels 4 and 5 are
of minimum size.
Figure~\ref{fig12} shows a minimum size skew army to reach level 6.

\begin{figure}[htb]
\centering
\epsfig{file=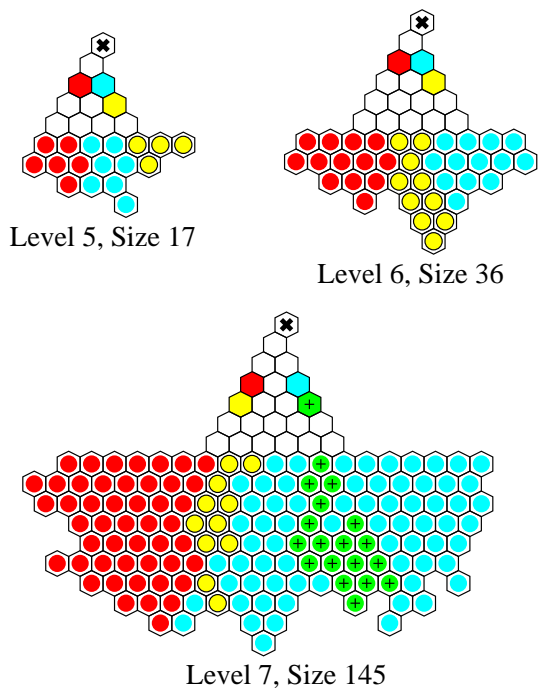}
\caption{Hexagonal armies to reach levels 5--7.  The first two are minimum size, the last can possibly be reduced by one man.}
\label{fig11}
\end{figure}

Finally, Figure~\ref{fig11} shows minimum size hexagonal armies.
The level 7 army is divided into five regiments.
The two blue regiments are separated by a green regiment,
and either blue regiment is capable of reaching
the blue target cell.
This is the only minimum army that was converted directly
from an IP solution.
For a more symmetrical army, which involves more men
(and hence is not minimum size), see \cite{Duncan}.
See \cite{Bellweb} for diagrams showing how
each army can advance to its specified level.

\subsection*{\normalsize 5.2 The hexagonal army to level 7}

Advancing \textit{any size} hexagonal army to level 7 is the
most difficult solitaire army problem in this paper,
and the reader should try it by hand.
In our experience, it cannot be done without
significant advance planning.

The determination of the minimum size hexagonal
army to level 7 is a difficult IP model to solve.
A direct run of the IP model presented above ran for ten hours,
and returned the 145-man solution shown in Figure~\ref{fig11},
but only gave a lower bound of 143 on the minimum army size.

As before, let $n_i$ be the number of men with exponent $i$.
Consider what the sequence $n_i$ must look like for a minimum size army.
First, it is impossible to have $n_i=0$ anywhere in the middle of the
army, with $n_{i-1}>0$ and $n_{i+1}>0$.
For if this were the case, then there would have to be some jump in the solution
that places a man into a cell with exponent $i$.
After this jump, the army is smaller but can still reach the same level,
hence the original army was not of minimum size.
For similar reasons, it must be that, for the maximum exponent $E$,
$n_{E-1}\ge n_E$.
In addition, the value of the slack cannot be too small,
by Proposition~1.
Table~\ref{table8} shows the possible values of the sequence $n_i$
for a minimum size army of 143 or 144 men.

\begin{table}[htbp]
\begin{center} 
\begin{tabular}{ | c | c | c | c | c | c | c | l | l |}
\hline	
Case & Men & $n_{13}$ & $n_{14}$ & $n_{15}$ & $n_{16}$ & $n_{17}$ & Slack & IP run time \\
\hline
\hline
1 & 143 & 25 & 23 & 8 & 0 & 0 & $0$ & 40 sec\\
\hline
2 & 144 & 25 & 22 & 9 & 1 & 0 & $0$ & $> 10$ hrs \\
\hline
3 & 144 & 24 & 24 & 9 & 0 & 0 & $0$ & $> 10$ hrs \\
\hline
4 & 144 & 25 & 23 & 7 & 1 & 1 & $0$ & 4 min \\
\hline
5 & 144 & 25 & 24 & 8 & 0 & 0 & $\sigma^{14}$ & 30 min \\
\hline
6 & 144 & 25 & 23 & 9 & 0 & 0 & $\sigma^{15}$ & 50 min \\
\hline
Figure~\ref{fig11}c & 145 & 25 & 21 & 10 & 2 & 0 & $0$ & \\
\hline
\end{tabular}
\caption{Number of men of each exponent,
and the slack, for the hexagonal army to reach level 7.
For all cases $n_i=3i-14$ for $7\le i\le 12$.
}
\label{table8}
\end{center} 
\end{table}

We can run the IP model on each of these cases separately by adding
additional constraints.
Suppose $M_i$ is the subset of cells with exponent $i$.
The added constraints are of the form
$$\sum_{m\in M_i} \mbox{STA}_m = n_i.$$

Several other ideas can be used to speed up the IP model.
First, if the slack is zero, there is no reason to include
any jump that loses weight, which reduces the jump
set $N$ by a factor of 3.
Instead of targeting a one-man finish,
we set $\mbox{FIN}_m$ as the 3-man configuration in
Figure~\ref{fig7}a.
This breaks the left-right symmetry of the original
problem and the IP solver does not have to search through
mirror symmetric board positions.

Cases 5 and 6 solve more quickly with the additional constraint
that there is only one jump that loses exactly $\sigma^{j}$,
where $j=14$ or $15$.
If $N_j$ is the subset of jumps $N$ that lose exactly $\sigma^j$,
then these constraints are of the form
$$\sum_{n\in N_j} J_n = 1.$$

With all of these improvements,
some of the cases in Table~\ref{table8} run relatively quickly,
returning ``global search complete---integer infeasible".
But two in particular are problematic,
and do not finish a global search in ten hours of run time on the NEOS solvers.
For this reason, although we believe a 144-man solution to be unlikely,
it cannot be ruled out.

\vskip 30pt
\section*{\normalsize 6. Summary and future work}

We have determined the minimum size army to reach
all levels for all five army types, with the exception
of the hexagonal army to level 7.
These results are summarized in Table~\ref{table7}.
Note that three of the columns in
Table~\ref{table7} are sequences in the
On-Line Encyclopedia of Integer Sequences \cite{OEIS}.

\begin{table}[htbp]
\begin{center} 
\begin{tabular}{ | c | c | c | c | c | c | }
\hline
 & \multicolumn{5}{|c|}{Minimum army size to reach level $n$} \\		
Level \# ($n$) & Conway's & Skew & Pablito's & Hexagonal & Diagonal  \\
\hline
OEIS id: & A014225 & ~ & ~ & A014227 & A125730 \\
\hline
\hline
1 & 2 & 2 & 2 & 2 & $2=F_3$ \\
\hline
2 & 4 & 3 & 3 & 3 & $3=F_4$ \\
\hline
3 & 8 & 5 & 5 & 5 & $5=F_5$ \\
\hline
4 & 20 & 9 & 9 & 9 & $8=F_6$ \\
\hline
5 & Impossible & 19 & 19 & 17 & $13=F_7$ \\
\hline
6 & Impossible & 46 & 53 & 36 & $23>F_8$ \\
\hline
7 & \multicolumn{3}{|c|}{Impossible} & 144 or 145 & $46$ \\
\hline
8 & \multicolumn{4}{|c|}{Impossible} & $123$ \\
\hline
9 & \multicolumn{5}{|c|}{All Impossible} \\
\hline
\end{tabular}
\caption{The minimum army size required to reach level $n$.} 
\label{table7}
\end{center} 
\end{table}

We have presented a number of techniques for
determining minimum size armies.
We have found useful the improved pagoda function technique,
as well as dividing the army into regiments whose target is
to reach a particular cell, short of the goal.
The IP model is a powerful technique for obtaining tight
lower bounds on army size.

The armies presented in this paper have minimum size,
but their shape is not unique.
It is known that, for Conway's army,
there are exactly four different configurations of 20
men that can reach level four (plus their reflections) \cite{Aigner}\footnote{\cite{WinningWays}
and \cite{Beasley} each state that up to reflection there are only two solutions.
However, the solutions given by each are fundamentally different
(not reflections of one another)!}.
It would be interesting to count the number of minimum size
armies for each entry in Table~\ref{table7}.

We could also consider how ``quickly" an army can advance.
All 20-man Conway's armies that reach level four do so in
exactly 19 jumps, because we start with 20 men, end with one,
and one man is lost per jump.
However, we can define a \textbf{move} as one \textit{or more}
jumps by the same man,
and it is not hard to verify that only ten moves are
needed to advance the fastest 20-man army.

It is possible to advance the 53-man Pablito's army shown
in Figure~\ref{fig10}c
using only 26~moves \cite{Bellweb}
(a significant improvement over the 73-man army
advancing in 43 moves given in \cite{Pablo1}).
It is not known whether this is the smallest number of moves
for this army (or any 53-man army) to reach level 6.
These move minimization problems appear very difficult for
the largest armies in this paper.

\vskip 15pt
\section*{\normalsize Acknowledgments}
We dedicate this paper to John H. Conway in
honor of his upcoming 70th birthday.
We also thank the anonymous referee for several suggestions
that improved the paper.

\end{document}